\documentclass{article}
\usepackage[T2A]{fontenc}
\usepackage[cp1251]{inputenc}
\usepackage[english]{babel}

\begin{document}
\begin{center} {\bf \Large
Stochastic  models for fully coupled systems of nonlinear parabolic equations }
\end{center} \footnote{Key words:
 stochastic differential equations, Markov chains
systems of semilinear parabolic equations, classical solutions of the Cauchy problem }
\begin{center} Belopolskaya Ya.
 \end{center}

\begin{abstract}

We construct a probabilistic representation of a system of fully
coupled parabolic equations arising as a model describing spatial segregation of interacting population species.  We derive a closed system of
stochastic equations such that its solution allows to obtain a
probabilistic representation of a weak solution of the Cauchy
problem for the PDE system. The corresponded stochastic system is
presented in the form of a  system of
stochastic equations  describing nonlinear Markov processes and their multiplicative functionals. 
\end{abstract}

\section{Introduction}

Systems of nonlinear  parabolic equations arise as models in
description of various physical, chemical and biological phenomena.
To mention some of them   recall such  phenomena as  chemotaxis   which is a biological
phenomenon describing the change of motion of a population densities
or of single particles in response (taxis) to an external chemical
stimulus spread in the environment  \cite{Per1}. Another example, are   cross-diffusion systems which describe 
segregation  processes in dynamics of interacting populations \cite{Jung}, \cite{SKT}
and in general a number of  diffusive conservation laws. 
Mathematical models of these conservation laws  are presented  as
systems of quasilinear parabolic equations  in both divergent 
and   non divergent form.

Actually, from the probabilistic point of view systems of nonlinear  parabolic equations can be divided into several large classes.Let us mention here two of them, namely, systems of type 1   \cite{LSU}
having diagonal entries  of second order terms  
 \begin{equation} \label{1.2}
{\bf u}^q_t=\frac 12\sum_{i,j=1}^dF^{ij}_u(x)\nabla^2_{ij} u^q+ \sum_{i=1}^da^i_u(x) \nabla_i u^q+\sum_{m=1}^{d_1} \sum_{i=1}^{d}[B_u]_i^{mq}(x)\nabla_i u^m+\sum_{m=1}^{d_1} 
c^{mq}_u(x)u^m,
\end{equation}
where coefficients $a_u, B_u, c_u$ and $F^{ij}_u=\sum_{k=1}^dA^{ki}_uA_u^{jk}, $   have the form $a_u(x)=a(x, u(x), \nabla u(x))$ 
and systems of type 2  \cite{Am}   with nondiagonal entries  of second order terms   such as
\begin{equation} \label{1.3}
{\bf u}^q_t=div G^q+ f^q(u), l=1,2,
\end{equation}
where $G^q=\alpha^q\nabla u^q+\beta^q\nabla (u^q)^2+\gamma^q\nabla [u^1u^2]$   and $\alpha_l,\beta_l,\gamma_l$ are constants.

The results concerning weak and classical solutions of the Cauchy problem for type 1 parabolic systems  could be found in the monograph
 \cite{LSU}   and a number of more recent papers. 
Probabilistic approach to investigation of different types of the  Cauchy problem  solution (namely, classical, generalized and viscosity)  for nonlinear parabolic systems of type 1  was  developed  for many years.    In particular, probabilistic representations for classical solutions of the Cauchy problem for nonlinear parabolic systems of this kind were constructed in  \cite{BD1} -\cite{BD2},   for generalized solutions in \cite{BW2} and for viscosity solutions in \cite{Bel2}.

Investigation of parabolic systems of  type 2 started  by Amann
\cite{Am}  was developed by many researches   especially in
connections with problems arising in applications.  These systems  are  sometimes  called  cross-diffusion systems  nevertheless there are only few papers devoted
  to  the investigation of their stochastic origin
 \cite{FM}, \cite{Gal}, \cite{Jung2} or the stochastic representation of their solutions
\cite{Bel3}.

 This paper is aimed to develop ideas  stated in   \cite{Bel2}, \cite{Bel3}.

We consider a particular case of the above system with cross diffusion used to model the
segregation of interacting species proposed in \cite{SKT}.
 The system under consideration has the form 
\begin{equation}\label{1.4}
 u^1_t=\Delta(  u^1 ( d_1+ d_{11}u^1+d_{12}u^2)) + 
u^1(\alpha_1- \alpha_{11}  u^1 -\alpha_{12}  u^2),\quad  u^1(0,x)=u^1_0(x),
\end{equation}
\begin{equation}\label{1.5}
u^2_t=\Delta (u^2 ( d_2+d_{21}u^1+d_{22}u^2)) +
u^2(\alpha_2-\alpha_{21} u^1 -\alpha_{22} u^2),\quad 
u^2(0,x)=u^2_0(x),
\end{equation}
where
$u^1$   and  $u^2$   are the densities of two competing  species. 

 Here and below we use notations $u_t=\frac{\partial u}{\partial t}$ for a partial derivative of
 $u(t,x)$ in $t\in[0,T]$ and $\nabla u=(\frac{\partial u}{\partial x_1},\dots,\frac{\partial u}{\partial x_d})$
 for  the gradient of $u$.

Note that a probabilistic approach to investigation of  the  parabolic equations includes there steps. The first is under assumption that there exists a unique solution
of the PDE problem (as a rule rather regular) one have to  derive a potential probabilistic representation in terms of some diffusion processes. The second step is to
derive a closed  stochastic system for these diffusion processes and  the third step is to investigate the stochastic problem and construct its solution with the required properties
which allow to check that it gives rise to the solution of the original PDE problem.

In our previous papers  this program was partly realized. 
The aim of this paper is to  derive a probabilistic
interpretation of a  generalized solution  of the
Cauchy problem for this system and to construct a closed stochastic system of stochastic equations.  

A probabilistic approach to construct a representation of a generalized solution to the
Cauchy problem for a scalar nonlinear parabolic equation based on the results of  the Kunita stochastic flows theory \cite{Ku}--\cite{Ku2}  was
developed  in 
\cite{BW1},  \cite{BW2}.

The main problem we meet  while extending this approach to systems with cross diffusion  is the fact that we cannot 
apply directly the Ito formula which is a main tool in construction of probabilistic representation of a generalized solution of  
the Cauchy problem  for scalar nonlinear parabolic equations \cite {BW1} and systems of nonlinear parabolic equations with diagonal principal part \cite{BW2}.

Nevertheless using the stochastic flow approach we succeed to derive the probabilistic representation of a generalized solution to
 (\ref{1.4}), (\ref{1.5}).

\section{ Probabilistic models for a fully coupled system of PDEs}
\setcounter{equation}{0}

\subsection{Dual PDEs}

To define a generalized solution of the Cauchy problem
(\ref{1.2}), (\ref{1.3}) we  have to introduce a number of functional
spaces.

$\bullet$   The set $C^k_b(R^d; R^{d_1})$   of $C^k$-smooth functions
$f: R^d\to R^{d_1}$  whose partial derivatives of order less or
equal then $k$ are bounded.

$\bullet$ The set $C^{1,k}_0([0,T]\times R^d, R^{d_1})$  of
continuous functions $u:[0,T]\times R^d\to R^{d_1}$ with compact
supports  with continuous  first partial derivatives in $t$ and partial
derivatives of order less or equal then $k$   in $x\in R^d$.

$\bullet$  The set $L^2=L^2(R^d;R)$  of real valued
functions with the norm
$\|u\|_{L^2}=\left(\int_{R^d}|u(x)|^2dx\right)^{\frac 12}$ which is  the Hilbert space with the
inner product $\langle u,v\rangle=\int_{R^d}u(x)v(x)dx$ for $u,v\in
L^2$.

$\bullet$    The Sobolev space   $ H^{k}=\{ u\in L^2: \nabla^ku \in
L^2\}$, where $\nabla^{k} u(x)$ is a generalized $k$-th order
derivative  with the norm
 $$\| u\|_{k}=\left(\int_{R^d}[|u(x)|^2+  \sum_{|\alpha|\le k}\|\nabla^\alpha u(x)\|^2]dx\right)^{\frac 12},$$
 where $|\alpha|= \sum\alpha_j$.

 $\bullet$  The Schwartz space ${\cal D}=C^\infty_0(R^d)$   and the
 space
 $ H^{k}_0= H^k_0$ which is a completion of $\cal D$ in the norm  $\|\cdot\|_k.$
 Given $ R^2$ -valued functions  $u, h$   with components in   ${\cal D}$ set
$$\langle\langle u,h\rangle\rangle  =\int_R u(x)\cdot h(x)dx,$$  where $ u\cdot h=u_1h_1 +u_2h_2.$

$\bullet$ ${\cal W}_{T}=\{u:u\in L^2((0,T);L^2(R^d))\cap
L^2((0,T);H^1_0(R^d))\}$ -- the set of   functions $t\mapsto
u(t,\cdot)$   from $(0,T)$ to   $L^2(R^d)$ with the norm
$$\|u\|_{\cal W}^2=\int_{R^d}\int_0^T[\|u(t)\|^2_{L^2}+\|\nabla u(t)\|^2_{L^2}]dt.$$

A couple of functions $  u^1,u^2\in {\cal W}_{T}$
is a weak (generalized) solution of the Cauchy problem
(\ref{1.4}), (\ref{1.5}) provided
 $  u^1, u^2 \in {\cal W}_T$    and equalities
\begin{equation}\label{2.1}
 \partial_s\int_{R^d} u^q(s,x)h(x)dx =\int_{R^d}u^q(s,x)[M^q_u
\Delta h(x)+c^q_{ u^1, u^2}(x)]h (x)dx
\end{equation}
hold for arbitrary  $h\in C_0^{1,\infty}([0,T]\times R^d)$   and $q=1,2$. Here
\begin{equation}\label{2.2}
M^q_u=d_q+d_{q1}u^1+d_{q2}u^2,\quad 
c^q_{u}=\alpha_q- \alpha_{q1}u^1-\alpha_{q2} u^2.\end{equation}

Along with this definition we need  one more  which gives a prompt concerning to what sort of diffusion processes we should be interested in.

We say that a couple of functions $  u^1,u^2\in {\cal W}_{T}$
is a weak (generalized) solution of the Cauchy problem
(\ref{1.4}), (\ref{1.5}) provided
 $  u^1, u^2 \in {\cal W}_T$    and equalities
\begin{equation}\label{2.1}
 \int_0^t\langle  u^q(s), \left\{ \partial_sh +\frac 12 [M^q_u]^2
\Delta h(s)+c^q_{ u^1, u^2}(s)h (s)\rangle
\right\}ds
\end{equation}$$=\langle  u^q(t),h(t)\rangle-\langle
 u^q(0), h(0)\rangle$$
hold for arbitrary  $h\in C_0^{1,\infty}([0,T]\times R^d)$   and $q=1,2$. Here
\begin{equation}\label{2.2}
M^q_u=\sqrt{2[d_q+d_{q1}u^1+d_{q2}u^2]},\quad 
c^q_{u}=\alpha_q- \alpha_{q1}u^1-\alpha_{q2} u^2.\end{equation}

We say that a couple  $( u^1,  u^2)\in {\cal W}^2_{s,T}$ is
a regular solution of the Cauchy problem (\ref{1.4})--(\ref{1.5})
provided $( u^1,  u^2)\in C^{1,2}([0,T]\times
R^d;R^2)\cap{\cal W}^2_T$ in addition to the above integrals
identities.

It should be noted that  similar to \cite{BRS1} we  can prove the following statement.

{\bf Lemma 2.1.} {\em  Let $u^q(t,x), q=1,2$ be a solution to (\ref{1.4})--(\ref{1.5}) such that  $\sup_{\theta\in (0,T)}|u^q(t,x)|<\infty$   and $u^q\in L^1([0,T]\times G) $  for each compact set  $G\in R^d$. Then   (\ref{2.1}) holds for a.a. $t\in [0,T]$ and  for every function $h\in C_b([0,T]\times R^d)\cap C_b^{1,2}((0,T)\times R^d)$
}

Proof.  It is enough to check this equality for any $h(t,x)=0$ for $\|x\|>R$ for some $R>0$. Let $\rho\in C^\infty_0((0,T))$.  
Since $u^q$ satisfies   (\ref{1.4})--(\ref{1.5}),
then applying integration by part formula we easily check that
$$\int_0^T\int_{R^d}u^q(t,x)[\partial_t(\rho h)+ {\cal M}^q(\rho h)](t,x)dxdt=0,$$
where
 \begin{equation}\label{2.3}{\cal M}^q h(t,x)=M^q_u
\Delta h(t)+c^q_{ u}(t)h.\end{equation}
Hence
$$-\int_0^T\frac{d}{dt}  \rho(t)\int_{R^d}u^q(t,x)h(t,x) dx\,dt=\int_0^T\rho(t)\int_{R^d}u^q(t,x)[\partial_t h(t,x)+{\cal L}h(t,x)]dx\,dt.$$
This yields that  the function  $t\mapsto \int_{R^d}u^q(t,x)h(t,x)dx$  has an absolutely continuous version and 
$$\frac{d}{dt} \int_{R^d}u^q(t,x)h(t,x)dx=\int_{R^d} u^q(t,x)[\partial_th(t,x)+ {\cal M}^q_u h(t,x)]dx.$$
Finally for some real number $C$
$$ \int_{R^d}u^q(t,x)h(t,x)dx=C+\int_0^t \int_{R^d}u^q(\theta,x)[\partial_\theta h(\theta,x)+{\cal M}^q_uh(\theta,x)dxd\theta \quad \mbox{ for a.a.}\, t\in [0,T].$$
Since by assumption $h(t,x)$ converges uniformly to $h(0,x)$ as $t\to 0$ we get
$$\lim_{t\to 0} \int_{R^d}u^q(t,x)h(0,x)dx\to \int_{R^d}u^q_0(x)h(0,x)dx$$
and thus $C=\int_{R^d}h(0,x)u^q_0(x)dx.$

Assume that there exists a unique $C^2$-regular generalized solution  $( u^1, u^2)$
to (\ref{1.4})-(\ref{1.5}), such that
 $\inf_x u^1(t,x)\ge 0, \inf_x u^2(t,x)\ge 0$.

To construct  stochastic processes associated with the Cauchy problem
  \begin{equation}\label{2.04}
\partial_sh^q +M^q_u(x)
\Delta h^q
+ c^q_{ u}(x)h^q =0,\quad h^q(T,x)=h^q(x),\quad 0\le s\le t \le T,
\end{equation}
we  consider
 a  stochastic differential equation (SDE)
    \begin{equation}\label{2.4}
d\xi^q(\theta)=  M^q_{u}(\xi^q(\theta))dw(\theta),   \quad
\xi^q(s)=y, \quad 0\le s\le \theta\le t\le T\end{equation}
and  a  linear SDE
\begin{equation}\label{2.5}
d\eta^q(\theta)=     c^q_{u}(\xi^q(\theta))\eta^q(\theta)d\theta,\quad
\eta^q(s)=1,\quad q=1,2.\end{equation}
Here
 $w(t)$   is a standard $R^d$-valued Wiener processes defined on
  a given probability space $(\Omega,{\cal F},P).$   
By the above assumption coefficients
  $M^q_u$ and $c^q_u$  in (\ref{2.4}), (\ref{2.5}) are  bounded Lipschitz continuous nonrandom functions. Hence, by standard reasoning we can
prove that there exists a Markov process
$\xi^q(t)$ satisfying (\ref{2.4}) and a multiplicative functional $\eta^q(t)$ (of the process $\xi^q(t)$) which satisfies  (\ref{2.5}). 
In  addition for $h\in C^2(R^d)\cap H^2(R^d)$ there exists a unique classical solution $h^q$ of the Cauchy problem (\ref{2.04})  that admits a probabilistic representation in the form
\begin{equation}\label{2.6}
h^q(s,y)=E[\eta^q(t)h^q(\xi^q_{s,y}(t))].
\end{equation}

Later we will see that to deal with the original Cauchy problem   (\ref{1.4})--(\ref{1.5}) we have to consider  a more complicated linear SDE of the form 
\begin{equation}\label{2.7}
d\eta^q(\theta)=   \tilde c^q_{u}(\xi(\theta))\eta^q(\theta)d\theta
+  C^q_{u}(\xi(\theta))\eta^q(\theta)dw(\theta),\quad
\eta^q(s)=1,\quad q=1,2,\end{equation}
with coefficients  $\tilde c^q_u=c^q_u+\|\nabla M_u\|^2$   and  $C_u=-\nabla M_u.$   

 Along with the process
$\xi^q(\theta)$ we will need the process $\hat\xi^q(\theta)=\xi^q(t-\theta)$. One can easily see that $\xi^q(t)=\hat\xi^q(0)$ and $\hat\xi^q(t)=\xi^q(0)$.

The  process $\xi^q(\theta)$   satisfying  (\ref{2.4}) and its time reversed $\hat\xi^q(\theta)=\xi^q(t-\theta)$   play an
important role  in construction of probabilistic representations for generalized 
 solutions of    (\ref{1.4})--(\ref{1.5}). 
 Similar to \cite{Bel3}  based
on the results of  \cite{Ku1},\cite{Ku2} we can prove   the following
statement.

{\bf Theorem 2.2} {\em  Given strictly positive regular (weak) solution $u^1,u^2\in {\cal W}_{0,T}\cap C^{1,2}$
 of the Cauchy problem  (\ref{1.4})--(\ref{1.5})
 let $\xi^q(\theta), \eta^q(\theta),\,q=1,2$  satisfy   (\ref{2.4}),  (\ref{2.7})
  while $\hat\xi^q(\theta)$  and  $\hat\eta^q(\theta)=\eta^q(t)\circ\psi^q_{\theta,t}$  be the corresponding time reversal processes.
  
   Then functions
\begin{equation}\label{2.8}
 u^q(s,x)=E_{s,x}\left[\hat\eta^q(t) u^q_0(\hat \xi^q(t))\right],\quad
 q=1,2,\end{equation}
 satisfy the integral identities
 \begin{equation}\label{2.9}
\langle u^q(t), h\rangle = \langle u^q(0), h\rangle -\int_0^t\langle [\Delta [ u^q(\theta)   M^q_u(\theta)] + 
c^q_u(x)u^q(\theta)]h\rangle \}d\theta
\end{equation}}

Note that the system  (\ref{2.4}) ,(\ref{2.7}) ,(\ref{2.8})  is not closed  since coefficients of SDEs governed 
 the processes $\hat\xi^q(\theta)$ and $\hat\eta^q(\theta)$ depend on $u^q$   and $ \nabla u^q$. Hence to make this system closed we need 
 probabilistic representations both for  functions $u^q$ and their gradients $\nabla u^q$  as well.

\subsection{Stochastic flows}

We start with  introducing  some required  functional spaces.

Given a Wiener process $w(\theta)$ denote by $\tilde w(\theta)=w(t-\theta)-w(t)$ for a fixed $t>\theta$ and given a stochastic process $\xi^q(\theta), s\le\theta\le t$  denote by  $
\hat\xi^q(\theta)=\xi^q( t-\theta)$  its time reversal process.    Let $
\phi^q_{0,\theta}:R^d\to R^d $,   $\psi^q_{0,\theta}:R\to R^d$ be
stochastic flows generated by  $ \xi^q(\theta)$  and $\hat\xi^q(\theta)$
that is $\xi^q_{0,y}(\theta)=\phi^q_{0,\theta}(y),\quad
\hat\xi^q_{0,x}(\theta)= \psi^q_{0,\theta}(x)$,  that is $\phi^q_{0,t}(\psi^q_{0,t}(x))=x$.

Given  bounded functions  $u^q(\theta,x)\in R^2, q=1,2,$  defined on $[0,T]\times R^d$ and
differentiable in $x$
we consider a stochastic equation  (\ref{2.4}).
 One can check (see \cite{Prot}, \cite{BD}) that the time reversal process  $\hat\xi_q^{0,x}(\theta)$ satisfies the equation
 \begin{equation}\label{2.10}
d\hat\xi^q_{0,x}(\theta)= [ M^q_u\nabla M^q_u](\hat\xi^q_{0,x}(\theta))d\theta+M^q_u(\hat\xi^q_{0,x}(\theta))d \tilde w(\theta ),\quad \hat \xi^q_{0,x}(0)=x.
\end{equation} where,  $0\le
\theta\le t$.

 In addition under the above assumptions  the mappings
  $\phi^q_{0,\theta}:x\mapsto \xi^q_{0,\theta}(x)$ and $\psi^q_{0,\theta}:y\mapsto \xi^q_{0,\theta}(y)$ are differentiable.   
  The process  ${\bf J}^q(t)=\nabla\xi^q(t)$ satisfies the SDE
   \begin{equation}\label{2.11}
d{\bf J}^q(\theta)= \nabla M^q_u(\theta) {\bf J}^q(\theta)dw(\theta),\quad {\bf J}^q(s)=I,
 \end{equation} 
 where $I$ is the identity matrix.
  Denote by
 $ J^q_{0,t}(\omega)={\rm det}{\bf  J}^q_{0,t}(\omega)$    and note that
  $$J^q_{0,t}(\omega)>0\quad\mbox{
and   }\quad  J^q_{s,s}(\omega))=1.$$
Moreover  $ J^q_{0,t}$ satisfies the SDE
 \begin{equation}\label{2.12}
dJ^q(\theta)= J^q(\theta) Tr [{\bf J}^q(\theta)d{\bf J}^q(\theta)],\quad  J^q(s)=1.
 \end{equation} 
The relation (\ref{2.12}) can be easily deduced from a known  formula for a differential of a determinant of a matrix  and  polylinearity
of the function $det{\bf J}^q$ w.r.t rows of the matrix ${\bf J}^q$. Note that for a linear function $f(x)$ we have $df(\xi(t))=\nabla f\cdot  d\xi(t)$ 
that is  there are no correction terms in   the Ito formula.

Along with SDE  (\ref{2.10}) we will need below an alternative SDE for the process $\hat \xi^q(\theta)$. To derive the required  SDE we  apply the Ito-Wentzell formula  
 to
the composition $\phi^q_{0,t}\circ\psi^q_{0,t}$.  Let
$$d\psi^q_{0,x}(\theta)=\lambda_u^q(\psi^q_{0,x}(\theta))d\theta+\Lambda^q_u(\psi^q_{0,x}(\theta))dw(\theta),\psi^q_{0,x}(0)=x,$$
where  coefficients $\Lambda^q_u$ and $\lambda^q_u$ should be defined to ensure that 
$\psi^q_{0,t}(x)$ is time reversal to $\phi^q_{0,t}(\kappa)$.  Since in this case   $\phi^q_{0,t}\circ\psi^q_{0,t}(x)=x$ then  we get
$$d[\phi^q_{0,t}\circ\psi^q_{0,t}(x)]=M^q_u(\psi^q_{0,\kappa}(t))dw(t)+[\nabla \phi^q_{0,t}](\psi^q_{0,t}(x))\Lambda^q_u(\psi^u_{0,t}(x))dw(t)$$
$$+[\nabla \phi^q_{0,t}](\psi^q_{0,t}(x))\lambda^q_u(\psi^q_{0,t}(x))dt+\nabla M^q_u(\psi^q_{0,t}(\phi^q_{0,t}(x)))\Lambda^q_u(\psi^q_{0,t}(x))dt=0.$$
Hence
$$\Lambda^q_u(\psi^q_{0,t}(x))=- [\nabla\phi^q_{0,t}(\psi^q_{0,t}(x))]^{-1}M^q_u(\psi^q_{0,t}(x)), $$$$\quad  \lambda^q_u(\psi^q_{0,t}(x))=
  [\nabla\phi^q_{0,t}(\psi^q_{0,t}(x))]^{-1}M^q_u(\psi^q_{0,t}(x))\nabla M^q_u(\phi^q_{0,t}(\kappa)).$$
Thus the process $\hat\xi^q(t)\equiv \hat\xi^q_{0,x}(t)=\psi^q_{0,t}(x)$ satisfies the SDE
 \begin{equation}\label{2.13}
 d\psi^q_{0,\theta}(x)=[\nabla\phi^q_{0,\theta}(x)]^{-1}M^q_u(x)\nabla M^q_u(x)d\theta- [\nabla\phi^q_{0,\theta}(x)]^{-1}M^q_u(x)dw(\theta),\quad \hat\xi^q(0)=x.\end{equation}
As a result  the following assertion holds.

{\bf Theorem 2.3} {\em Let $\phi^q_{s,t}(\kappa)$ be a solution of  
(\ref{2.9})
 with $C^k$ -smooth coefficients for  $(k\ge 3)$.
Then the inverse flow  $[\phi^q_{s,t}]^{-1}=\psi^q_{s,t}$  satisfies
 (\ref{2.13}) and is a diffeomorphism.
}
 
Given a generalized function  $u^q_0$ we define a composition
 of $u^q_0$  with a stochastic flow   $\psi^q_{s,t}(\omega)$   as a random variable valued in
   the space ${\cal D}'$ dual to ${\cal D}$. Note that given
   $h\in{\cal D}$  a product  $h\circ\phi^q_{s,t}(\omega)  J^q_{s,t}(\omega)$ belongs to   $\cal D$. Set
\begin{equation}\label{2.14}
\langle u^q_0,T^q_{s,t}h(\omega)\rangle=\langle u^q_0, h\circ\phi^q_{s,t}(\omega)
 J^q_{s,t}(\omega)\rangle,\quad h\in{\cal D}.
\end{equation}
One can check that (\ref{2.14}) defines a linear functional on
${\cal D}$ which we denote by  $u^q_0\circ \psi^u_{s,t}$. If $m$ can be
represented as    $u^q_0=u^q_0(x)dx,$ where $u^q_0(x)$ is a continuous function
then  $u^q_0\circ\psi^u_{s,t}$ is just a composition of $u_0$ with
$\psi^u_{s,t}$ due to an equality
\begin{equation}\label{2.15}
\int_{R^d}u^q_0(\psi^q_{0,t}(x,\omega))h(x)dx=
\int_{R^d}u^q_0(y)h(\phi^q_{0,t}(y,\omega))
 J^q(t))dy,\quad h\in{\cal D},
\end{equation}
resulting from the changing variable formula under the integral sign
$$
\langle u^q_0\circ \psi^q_{0,t}(\cdot),h\rangle=\int_{R^d} u^q_0(
\hat\xi^q_{0,x}(t))h(x)dx=\int_{R^d}u^q_0(y)h(\xi^q_{s,y}(t))
J^q(t)d\kappa.$$
By similar arguments we conclude that  when $ u^q(t)\circ  \psi^q_{0,t}(x)=\hat\eta^q_{0,x}(t)u_0^q(\hat\xi^q_{0,x}(t))$, 
$$
\langle  u^q(t) \circ \psi^q_{s,t}(\cdot),h\rangle=\int_{R^d}\hat\eta^q_{s,x}(t)u_0^q(\hat\xi^q_{s,x}(t))h(x) dx$$
$$=\int_{R^d}u_0^q(y) \eta^q_{s,y}(t)h(\xi^q_{s,y}(t))
 J^q(t))dy.$$

Let   
\begin{equation}\label{2.16}\tilde{\cal  M}^q_u(x)=\frac 12 [M^q_u]^2(x)\Delta,\quad {\cal M}^q_u=\tilde{\cal M}^q_u+ c^q_u\end{equation}
and 
 $ {\cal M}^q_u$ be dual operators to  ${\cal L}^q_u$  in $L^2(R^d)$,  that is
 $$\langle  {\cal L}^q_vu^q,
h\rangle=\langle u_0^q,{\cal M}^q_uh\rangle,\quad \langle  {\cal L}^q_u u^q,
h\rangle=\langle u^q,{\cal M}^q_uh\rangle.$$

By the above considerations we note that given a function $h\in {\cal D}$  one can consider  a  process  $\gamma^q(t)$  defined by  
\begin{equation}\label{2.17}
\gamma^q(t)= \eta^q(t)h(\xi^q_{0,y}(t)) J^q(t)
\end{equation}
and use $\gamma^q(t)$ as stochastic test functions.
Actually,   since    $\langle  \hat\eta^q(t)u^q_0\circ\psi^q_{0,t},h\rangle= \langle
u^q_0, \gamma^q(t)\rangle$  we can state the
following assertion.

{\bf Lemma 2.4.}{\em  Let coefficients $\tilde c^q_u$ and  $ C^q_u$   have the form
 \begin{equation}\label{2.18}
 \tilde c^q_u(\xi^q(\theta))=
c^q_u(\xi^q(\theta))-\langle \nabla M^q_u(\xi^q(\theta)),\nabla M^q_u(\xi^q(\theta))\rangle,\quad  C^q_u(\xi^q(\theta))=-\nabla M^q_u(\xi^q(\theta)).
\end{equation} 
Then  the processes $\gamma^{q}(\theta)=\eta^{q}(\theta)h(\xi^q_{0,y}(\theta)) J^q(\theta), q=1,2,$  have  stochastic differentials of the form
\begin{equation}\label{2.19}
d\gamma^{q}_y(\theta)=\left[  \frac 12[ M^q_u]^2\Delta h+c^q_uh\right](\xi^q_{0,y}(\theta))\eta^q(\theta)J^q(\theta)d\theta\end{equation} 
$$+[M^q_u\nabla h(\xi^q_{0,y}(\theta))  \eta^q(\theta)J^q(\theta)dw(\theta).
$$
}
Proof.
We apply the Ito formula to evaluate $d\gamma^q(t)$
$$
d\gamma^{q}(\theta)=d[\eta^{q}(\theta)h(  \xi^q(\theta)) J(\theta) ]=  d[\eta^{q}(\theta)]h( \xi^q(\theta))J^q(\theta) +\eta^{q}(\theta)d[h(  \xi^q(\theta))] J^q(\theta)$$$$
 +\eta^{q}(\theta)h( \xi^q(\theta)) dJ^q(\theta) 
+d[\eta^{q}(\theta)]d[h( \xi^q(\theta)) ]J^q(\theta) $$
$$+\eta^{q}(\theta)d[h( \xi^q(\theta))]d J^q(\theta) +d[\eta^{q}(\theta)]h( \xi^q(\theta)) dJ^q(\theta) .$$
Taking into account the expressions for $d\xi^q(t), d\eta^q(t)$ and $dJ^q(t)$   from   (\ref{2.4}),  (\ref{2.7})  and (\ref{2.12})  we deduce
$$
d\gamma^{q}(\theta)=\{ \tilde c^q_uh+\frac 12 [M^q_u]^2\Delta h+ \langle C^q_u,M^q\nabla h\rangle
+\langle \nabla M^q_u, M^q_u\nabla h\rangle\}(\xi(\theta))
\eta^q(\theta)J^q(\theta)d\theta$$$$+ \langle C^q_u, h\nabla M^q_u \rangle](\xi^q(\theta))
\eta^q(\theta)J^q(\theta)d\theta+\{C^q_uh +\nabla [M^q_u h] \}(\xi^q(\theta))\eta^q(\theta)J^q(\theta)dw(\theta).$$
Setting $$C^q_u(\xi(\theta))=-\nabla M^q_u(\xi(\theta))$$ and 
$$\tilde c^q_u(\xi(\theta))=
c^q_u(\xi(\theta))+\langle \nabla M^q_u(\xi(\theta)),\nabla M^q_u(\xi(\theta))\rangle $$   we get 
$$
d\gamma^{q}(\theta)=\left[  c^q_u(\xi^q(\theta))h(\xi(\theta))+\frac 12 [M^q_u]^2(\xi^q(\theta))\Delta h(\xi^q(\theta))\right]\eta^q(\theta)J^q(\theta)d\theta$$
$$+M^q_u(\xi^q(\theta))\nabla h(\xi^q(\theta)) \eta^q(\theta)J^q(\theta)dw(\theta).
$$

 To conclude this section  remind (see \cite{Ku1}) that given   functions  $m^q\in H^{1} $   we can define random variables  $\langle \eta^q(t)m^q\circ \psi^q_{0,t}, h\rangle$
 and prove that it has finite moments for $m^q\in H^{1}$ and $h\in H^{1}\cap C^\infty_0={\cal R}$.  Hence $\langle K^q_{0,t}, h\rangle=E\langle \eta^q(t)m^q\circ \psi^q_{0,t}, h\rangle$ is a continuous linear functional on $\cal R$
 that yields $K^q_{s,t}\in H^{1}$   and $K^q_{s,t}=E[ \eta^q(t)m^q\circ \psi^q_{s,t}]$  is called a generalized expectation of $g^q(t)=\eta^q(t)m^q\circ \psi^q_{s,t}$.
 Denote by $U^q(t)m^q=E[\eta^q(t)m^q\circ\psi_{0,t}]$. It is a linear map from $H^1$ into itself, possessing the semigroup property $U^q(t)U^q(s)=U^q(s+t)$  for any $s,t>0$.
 IOne can easily verify this fact  using stochastic flow properties.
 $$\langle U^q(t)U^q(s)m^q,h\rangle=E[\langle U^q(s)m^q, \eta^q(t)h\circ \phi^q_{s,t} J^q_{s,s+t}\rangle]$$$$=E[\langle m^q, \eta^q(s)h\circ\phi^q_{0,s}J_{0,s}\rangle |_{G^q=\eta^q(s+t)h\circ \phi^q_{s,s+t} \hat J^q_{s,s+t}}]= E[\langle m^q, (G^q\circ\phi^v_{s,s+t}\circ \phi^q_{s,0})(  J^q_{s,s+t} J^q_{0,s})\rangle ]$$
 $$= E[\langle m^q, G^q\circ\phi^q_{s+t,0}\hat J^q_{0,s+t}\rangle ]=\langle U^q(t+s)m^q,h\rangle.$$
 
 Due to evolution properties of multiplicative operator functionals we can check as well that   families  
 $U^q_{0,t}m^q=E[\eta^q(t)m^q(t)\circ\psi_{0,t}],$    have the  evolution property.

\section{Stochastic counterpart of the Cauchy problem for a system with cross-diffusion}
\setcounter{equation}{0}

Note that to obtain a closed system of stochastic equations which can be treated as a stochastic counterpart of the system 
\begin{equation}\label{3.1}
\frac{\partial u^q}{\partial t}=\Delta[u^q[d_q+d_{q1}u^1+d_{q2}u^2]]+c^q_uu^q,\quad u^q(0,x)=u_{0q}(x),\quad q=1,2.
\end{equation} 
it is not enough to construct a stochastic representation of the solution to (\ref{3.1}) itself since coefficients of SDEs for $\hat\xi^q(\theta)$ and 
$\hat\eta^q(\theta)$ derived in the previous section depend on $\nabla u^q$. Hence we need   stochastic representations for  spatial derivatives of $u^q(t,x)$.
In the previous section it was shown that we can associate   with (\ref{3.1})
  a system of stochastic equations
 \begin{equation}\label{3.2}
 d\xi^q_{0,y}(\theta)= M^q_u(\xi^q_{0,y}(\theta))dw(\theta ) ,\quad \xi^q_{0,y}(0)=y,
\end{equation} 
\begin{equation}\label{3.3}
   d\eta^q(\theta)=\tilde c^q_u(\xi^q_{0,y}(\theta))\eta^q(\theta)d\theta+ C^q_u(\xi^q_{0,y}(\theta))\eta^q(\theta)dw(\theta),\quad \eta^m(0)=1,\quad 
    \end{equation} 
 \begin{equation}\label{3.4}
 d\hat\xi^q_{0,x}(\theta)=  [M^q_u\nabla M^q_u](\hat\xi^q_{0,x}(\theta))d\theta +M^q_u(\hat\xi^q_{0,x}(\theta))d \tilde w(\theta ),\quad \hat\xi_{0,x}(0)=x,\end{equation} 
  \begin{equation}\label{3.5}
   d\hat\eta^q(\theta)=\tilde c^q_u(\hat\xi^q_{0,x}(\theta))\hat\eta^q(\theta)d\theta+C^q_u(\hat\xi^q_{0,x}(\theta))(\hat\eta^q(\theta),
   dw(\theta)),\
    \end{equation} 
    $$ \hat\eta^q(0)= 1, \quad q=1,2,$$
 \begin{equation}\label{3.6}
 u^q(t,x)=E[\hat\eta^q(t)u_0^q(\hat\xi_{0,x}(t))].
\end{equation} 
To make the system (\ref{3.4})-(\ref{3.6}) closed we need an extra relation for the function $v^q_i(t,x)=\nabla_i u^q(t,x)$, $i=1,\dots, d$,
since coefficients of (\ref{3.4}) and (\ref{3.5})  depend on $\nabla u^q$. 

To derive this relation we need some additional speculations based on results from [12]. 
By formal differentiation of the system
\begin{equation}\label{3.7}
\frac{\partial u^q}{\partial t} =\Delta [u^q(d_q+d_{q1}u^1+d_{q2}u^2)]+c^q_uu^q, \quad u^q(0,x)=u^{q}_0(x), \, q=1,2.
\end{equation} 
we get  a PDE   for  $v^q_i=\nabla_iu^q$
\begin{equation}\label{3.8}
\frac{\partial v^q_i}{\partial t} =\Delta\{v^q_i(d_q+d_{q1}u^1+d_{q2}u^2)+u^q(d_{q1}v^1+d_{q2}v^2)\}+ u^q\nabla_i c^q(u)+c^q(u)v^q_i, \quad 
\end{equation} 
$v^q_i(0,x)=\nabla_iu^q_0(x).$

In a similar way  from
\begin{equation}\label{3.9}
\frac{\partial h^q}{\partial \theta} + (d_q+d_{q1}u^1+d_{q2}u^2)\Delta h^q+c^q(u)h^q=0, \quad h^q(t,y)=h(y), \, 
\end{equation} 
 we get a PDE  for  $g^q_i=\nabla_i h^q$
\begin{equation}\label{3.10}
\frac{\partial g^q_i}{\partial \theta} + ((d_q+d_{q1}u^1+d_{q2}u^2))\Delta  g^q_i +(d_{q1}v^1_i+d_{q2}v^2_i)div  g^q+\nabla_ic^q(u)h^q+c^q(u)g^q_i=0, \quad \end{equation} 
$g^q_i(0,y)=\nabla_i h^q(y).
$

In addition note that we can construct  a stochastic representation of the solution to  (\ref{3.9})-(\ref{3.10})  in the form
$$G^q(\theta,y)=E[\beta^q(t) G_0(\xi^q_{\theta,y}(t))], $$
where  $G(t,y)=\pmatrix{h^q(t,y)\cr \nabla h^q(t,y)}$ and stochastic processes $\xi^q(\tau)$ and $\beta^q_{ik}(\tau)$   satisfy SDEs
$$d\xi^q(\tau)=\sqrt{2[d_q+d_{q1}u^1(t,\xi^q(\tau))+d_{q2}u^2(t,\xi^q(\tau))]}dw(\tau),\quad \xi^q(\theta)=y, 0\le\theta\le\tau \le t,$$
$$d\beta^q(\tau)=n^q_u(\xi(t))\beta^q(\tau)d\tau+ N^{q}_u(\xi(\tau))\beta^q(\tau)dw(\tau).$$
Here  
$$\beta^q(\tau)=\pmatrix{beta^q(\tau)&0\cr \nabla \beta^q(\tau)&\beta^q(\tau)},\quad     n^q_{u}=\pmatrix{c_u^q&0\cr \nabla c_u^q&c^q_u},\quad N^m_{u}=\pmatrix{0&0\cr 0&\frac{[d_{q1}v^1+d_{q2}v^2] \delta}{\sqrt{2(d_q+d_{q1}u^1+d_{q2}u^2)}}}$$
 where $\delta$ is the Kronecker delta.  Thus  for $G_0(y)=G^q(0,y)=\pmatrix{h^q(y)\cr\nabla h^q(y)}$ we get
$$G^q(\theta,y)=E\left[\pmatrix{\beta^q(t)&0\cr \nabla\eta^q(t)&\beta^q(t)}\pmatrix{h(\xi^q_{\theta,y}(t))\cr\nabla h(\xi^q_{\theta,y}(t))}\right]$$$$=
\pmatrix{ E[\beta^q(t)h(\xi^q_{\theta,y}(t))]\cr E[ \nabla\beta^q(t)h(\xi^q_{\theta,y}(t))+\beta^q(t)\nabla h(\xi^q_{\theta,y}(t))]}.$$

To deduce the stochastic representation  for the function  $g^q_j=\nabla_j h^q$ we note  that  given  the PDE system  (\ref{3.9})-(\ref{3.10}) 
we can derive its stochastic representation as follows. 
Let  us  rewrite the system  (\ref{3.7}),(\ref{3.8})   in the form 

\begin{equation}\label{3.11}
\frac{\partial }{\partial t}\pmatrix{u^q\cr v^q}= {\cal Z}^q\pmatrix{u^q\cr v^q},   \quad q=1,2.
\end{equation} 
where $$ {\cal Z}^q \pmatrix{u^q\cr v^q}=\Delta \left[\pmatrix{d_q+d_{q1}u^1+d_{q2}u^2&0\cr d_{q1}v^1+d_{q2}v^2&d_q+d_{q1}u^1+d_{q2}u^2}\pmatrix{u^q\cr v^q
}\right]+\pmatrix{c^q_{11}&0\cr c^q_{21}&c^q_{22}}\pmatrix{u^q\cr v^q}.$$

Consider as well a dual system derived from  (\ref{3.11}) as follows.  Integrate over      $R^d$ a  product of  (\ref{3.11})   and   vector test functions
$(h, g)^*$, where $g_j=\nabla_j h,\, j=1,\dots,d$. As a result we obtain a system of the form
\begin{equation}\label{3.12}
\left\langle\left\langle\pmatrix{u^q\cr v^q}\left[\frac{\partial}{\partial t}\pmatrix{h^q\cr g^q}+ {\cal Q}^q \pmatrix{h^q\cr g^q}\right]\right\rangle\right\rangle=0,
\end{equation} 
where  $${\cal Q}^q \pmatrix{h\cr g}=\pmatrix{d_q+d_{q1}u^1+d_{q2}u^2&0\cr d_{q1}v^1+d_{q2}v^2&d_q+d_{q1}u^1+d_{q2}u^2}
\Delta \pmatrix{h^q\cr g^q}+
 \pmatrix{c^q_{11}&0\cr c^q_{21}&c^q_{22}}\pmatrix{h^q\cr g^q}.$$
Here and below we denote by 
$$\left\langle\left\langle\pmatrix{u\cr v_i}\pmatrix{h\cr g_i}\right\rangle\right\rangle=\pmatrix{\int_{R^d}u(x)h(x)dx\cr\int_{R^d}v_i(x)g_i(x)dx}.$$
In the sequel we take into account the relation    $[d_{q1}v^1+d_{q2}v^2]\Delta h=[d_{q1}v^1+d_{q2}v^2] div g$ that allows to construct a proper stochastic representation of the backward Cauchy problem   
\begin{equation}\label{3.13}
\frac{\partial}{\partial \theta}\pmatrix{h^q\cr g^q}+ {\cal Q}^q \pmatrix{h^q\cr g^q}=0, \quad \pmatrix{h^q(t)\cr g^q(t)}=\pmatrix{h^q_0\cr g^q_0},\quad 0\le \theta\le t
\end{equation} 
based on results of [14].

 To this end along with   (\ref{3.2})     we consider a  stochastic equation of the form 
\begin{equation}\label{3.14}
d\beta^q(\theta)=[\tilde c^q]^*(\xi^q(\theta))\beta^q(\theta)d\theta+[ \tilde C^q]^*(\xi^q(\theta))(\beta^q(\theta),dw(\theta)),\quad 
\beta^q(s)=\beta^q\end{equation} 
with respect to  the two component process  $\beta^q(\theta)=\pmatrix{\beta^q_{1}(\theta)\cr\beta^q_{2}(\theta)}$    with coefficients  $\tilde c^q$ and $ C^q$    to be chosen below.
Let $\bar\beta^q(\theta)$ maps $\beta^q$ to $ \beta^q(\theta)$, that is 
$$\bar\beta^q(\theta)=\pmatrix{\bar\beta^q_{11}(\theta)&0\cr \bar\beta^q_{21}(\theta)&\bar\beta^q_{22}(\theta)}.$$

Define a  stochastic test function
\begin{equation}\label{3.15}
\kappa^q(\theta)=\pmatrix{\kappa^q_1(\theta)\cr\kappa^q_2(\theta)}=\pmatrix{\bar\beta^q_{11}(\theta)&0\cr \bar\beta^q_{12}(\theta)& \bar\beta^q_{22}(\theta)}
\pmatrix{h^q(\xi^q(\theta))\cr g^q(\xi^q(\theta))}J^q(\theta),
\end{equation} 
where $J^q(\theta)$  is a Jacobian of the stochastic transformation  $y\to\xi^q_{s,y}(\theta)$.  The  stochastic 
differential of the process $\kappa^q(\theta)$ has the form
$d\kappa^q(\theta)=\pmatrix{d\kappa^q_1(\theta)\cr d\kappa^q_2(\theta)}$ with
$$d\kappa^q_1(\theta)=[\tilde c^q_{11}h^q+\frac 12[M^q_u]^2\Delta h^q+\langle  C^q_{11},
[M^q_u\nabla h^q+\nabla M^q_uh^q] (\xi(\theta))\rangle \beta^q_{11}(\theta)J^q(\theta)]d\theta
$$
$$
+ \langle M^q_u\nabla h^q(\xi(\theta)),\nabla M^q_u\rangle \beta^q_{11}(\theta)J^q(\theta)d\theta+ \langle  N^q_1(\xi(\theta)), dw(\theta)\rangle,
$$
$$d\kappa^{iq}_2(\theta)=\left[[\tilde c^q_{21}h^q+ M^q_u\nabla M^q_u div g^q](\xi^q(\theta))\beta^{iq}_{21}(\theta)\right]J^q(\theta)d\theta$$$$+\left[\beta^q_{22}(\theta)[\tilde c_{22}g^q_i(\xi^q(\theta))+\frac 12[M^q_u]^2\Delta g^q_i(\xi^q(\theta))\right]J^q(\theta)d\theta$$
$$+ \{ C^q_{21}\beta^{iq}_{21}(\theta)[M^q_u\nabla h^q(\xi^q(\theta))+\nabla M^q_uh^q(\xi(\theta))]+ C^q_{22}\beta_{22}(\theta)[M^q_u\nabla g^q_i+
g^q_i\nabla M^q_u](\xi(\theta))$$
$$+\beta^{iq}_{21}(\theta)M^q_u \langle \nabla h^q,\nabla M^q_u\rangle(\xi^q(\theta))+\beta^q_{22}(\theta)M^q_u\langle \nabla g_i,\nabla M^q_u\rangle (\xi^q(\theta))\}J^q(\theta)d\theta$$$$+\langle [N_{21}(\xi^q(\theta))\beta^{iq}_{21}(\theta)+N^{iq}_{22}(\xi^q(\theta))\beta^q_{22}(\theta)], dw(\theta)\rangle J^q(\theta).
$$
  Let us specify coefficients $\tilde c^q$ and $C^q$. As it was done in the previous section we choose 
\begin{equation}\label{3.16}
 C^q_{11}=-\nabla M^q_u ,\quad  \tilde c^q_{11}=c^m_u+\|\nabla M^q_u\|^2.\end{equation}  
 Next we choose
\begin{equation}\label{3.17}
 C^q_{21}=-\nabla M^q_u,  C^m_{22}= \frac{ (d_{q1}v^1+d_{q2}v^2) \delta}{M^q_u} - \nabla M^q_u, \end{equation} 
 \begin{equation}\label{3.18}
  [\tilde c^q_{21}]_i=\nabla_i c^q_u+\|\nabla M^q_u\|^2,\quad \tilde c^q_{22}=
 c^q_u+\|\nabla M^q_u\|^2.
\end{equation} 
 We do not specify for the moment $N^q_1$ and $N^q_2$ since they do not take part in the probabilistic representation of $u^q$
  and $v^q$.

Next we proceed as in the previous section.  Denote by $\hat\beta^q(\theta)$ the time reversal of the process $\bar\beta^q(\theta)$.
To get a closed counterpart of the system (\ref{1.2}) in addition to theorem 1 we state the following assertion.

{\bf Theorem 3.1.}{\em Under assumptions of theorem 1  both the functions $u^q(t,x)$ admit stochastic representations (\ref{3.10}) and functions $v^q_j=\nabla_j u^q$
admit  stochastic representations
\begin{equation}\label{3.19}
\pmatrix{u^q(t,x)\cr\nabla_i u^q(t,x)}=E\left[ \pmatrix{\hat\beta^q_{11}(t)&0\cr \hat\beta^q_{21}(t)&\hat\beta^q_{22}(\theta)} \pmatrix{u_0^q(\hat\xi^q_{0,x}(t))\cr v^q_{i}(\hat\xi^q_{0,x}(t))}\right].
\end{equation} 
}

Proof.  To verify the last assertion of the theorem we  note that  we have the following matrix relations
\begin{equation}\label{3.20}
\left\langle \left\langle\int_0^t \pmatrix{u_0^q\cr v^q_{i0}} \pmatrix{d\kappa^q_1(\theta)\cr d\kappa_2^q(\theta)} \right\rangle\right\rangle  =\left\langle \left\langle \pmatrix{u_0^q\cr v^q_{i0}} \pmatrix{d\kappa^q_1(t)\cr d\kappa_2^q(t)}  \right\rangle\right\rangle\end{equation} $$-\left\langle\left\langle  \pmatrix{u_0^q\cr v^q_{i0}} \pmatrix{d\kappa^q_1(0)\cr d\kappa_2^q(0)}\right\rangle\right\rangle.$$
At the other hand from  (\ref{3.15})    we deduce
\begin{equation}\label{3.20}
E\left[\left\langle \left\langle\int_0^t \pmatrix{u_0^q\cr v^q_{i0}}   \pmatrix{d\kappa^1(\theta)\cr d\kappa^2(\theta)}  \right\rangle\right\rangle\right]\end{equation} 
$$=
E\left[\int_0^t\left\langle \left\langle  \pmatrix{u_0^q\cr v^q_{i0}}  d\left[  \pmatrix{\beta^m_{11}(\theta)&0\cr \beta^q_{21}(\theta)&\beta^q_{22}(\theta)}  \pmatrix{h^q(\xi^q_{0,\cdot}(\theta))\cr g^q(\xi^q_{0,\cdot}(\theta))}  J^q(\theta)\right]\right\rangle\right\rangle
 \right]
 $$
 $$
 =E\left[\int_0^t\left\langle \left\langle \pmatrix{u_0^q\cr v^q_{i0}}  \pmatrix{\beta^q_{11}(\theta)&0\cr \beta^q_{21}(\theta)&\beta^q_{22}(\theta)}  {\cal Q}^q
  \pmatrix{h^q(\xi^q_{0,\cdot}(\theta)\cr g^q(\xi^q_{0,\cdot}(\theta))}  J^q(\theta) \right\rangle \right\rangle  d\theta\right].
$$
By the change of variables $\xi_{0,y}(\theta)=x$ applying stochastic Fubini theorem we get 
\begin{equation}\label{3.21}
E\left[\left\langle \left\langle \int_0^t \pmatrix{u_0^q\cr v^q_{i0}} , \pmatrix{d\kappa^1(\theta) \cr d\kappa^2(\theta)}  \right\rangle \right\rangle \right]\end{equation} $$
=E\left[\int_0^t\left\langle \left\langle\pmatrix{\hat\beta^q_{11}(\theta)&0\cr \hat\beta^q_{21}(\theta)&\hat\beta^q_{22}(\theta)}  \pmatrix{u_0^q(\hat\xi^q_{0,\cdot}(\theta))\cr v^q_{i0}(\hat\xi^q_{0,\cdot}(\theta))} {\cal Q}^q\pmatrix{h^q\cr g^q} \right\rangle \right\rangle d\theta\right]
$$
$$=
\int_0^t\left\langle \left\langle E\left[\pmatrix{\hat\beta^q_{11}(\theta)&0\cr \hat\beta^q_{21}(\theta)&\hat\beta^q_{22}(\theta)} \pmatrix{u_0^q(\hat\xi^q_{0,\cdot}(\theta))\cr v^q_{i0}(\hat\xi^q_{0,\cdot}(\theta))} \right]{\cal Q}^m\pmatrix{h^q\cr g^q} \right\rangle \right\rangle d\theta$$
$$= \int_0^t\left\langle \left\langle{\cal Z}^qE\left[\pmatrix{\hat\beta^q_{11}(\theta)&0\cr \hat\beta^q_{21}(\theta)&\hat\beta^q_{22}(\theta)} \pmatrix{u_0^q(\hat\xi_{0,\cdot}(\theta))\cr v^q_{i0}(\hat\xi_{0,\cdot}(\theta))} \right]\pmatrix{h^q\cr g^q}  \right\rangle \right\rangle  d\theta.$$
Hence we derive that the functions $$\pmatrix{\lambda^q(t,x)\cr\nabla\lambda^q(t,x)} =E\left[\pmatrix{\hat\beta^q_{11}(\theta)&0\cr \hat\beta^q_{21}(\theta)&\hat\beta^q_{22}(\theta)}  \pmatrix{u_0^q(\hat\xi^q_{0,x}(\theta))\cr v^q_{i0}(\hat\xi_{0,x}(\theta))} \right]$$ satisfy integral identities
$$
\left\langle \left\langle\pmatrix{\lambda^q(t)\cr\nabla \lambda^q(t)} \pmatrix{h^q\cr g^q} \right\rangle \right\rangle - 
\left\langle \left\langle\pmatrix{\lambda^q(0)\cr\nabla \lambda^q(0)} \pmatrix{h^q\cr g^q} \right\rangle \right\rangle 
$$
$$
=\left\langle \left\langle{\cal Q}^q\pmatrix{\lambda^q(t)\cr\nabla \lambda^q(t)} \pmatrix{h^q\cr g^q} \right\rangle \right\rangle 
$$
which results due to the assumed uniqueness of a solution to (\ref{1.2}) that  
$$\pmatrix{\lambda^q(t,x)\cr\nabla \lambda^q(t,x)}=\pmatrix{u^q(t,x)\cr\nabla u^q(t,x)} $$
and hence
$$ \pmatrix{u^q(t,x)\cr\nabla u^q(t,x)}=E\left[\pmatrix{\hat\beta^q_{11}(t)&0\cr \hat\beta^q_{21}(t)&\hat\beta^q_{22}(t)}  \pmatrix{u_0^q(\hat\xi_{0,x}(t))\cr v^q_{i0}(\hat\xi^q_{0,x}(t))} \right].
$$
As a result  we deduce from the last equality that (\ref{3.14})   holds and in addition
$$\nabla u^q_i(t,x)=E[\hat\zeta_{21}(t)u_0^q(\hat\xi_{0,x}(t))+\hat\zeta^m_{22}(t) v^q_{i0}(\hat\xi^q_{0,x}(t))].$$

{\bf Remark.}
We have proved that under a priori assumption about  existence of a unique regular solution of the Cauchy problem (\ref{1.2}) there exists a stochastic representation of this solution and moreover we derive a closed system of stochastic equations that can be considered separately  without any reference to this a priori assumption. To be more precise  we have shown  that the system  (\ref{3.4}),  (\ref{3.14}) and  (\ref{3.18})  with coefficients given by  (\ref{3.15}) -- (\ref{3.17})  is a  closed stochastic system which can be considered independently of (\ref{1.2}). 

At the next step starting with the  system    (\ref{3.4}),  (\ref{3.14}) and  (\ref{3.19})  which can be rewritten in the form
$$
 d\hat\xi^q_{0,x}(\theta)=  [M^q_u\nabla M^q_u](\hat\xi^q_{0,x}(\theta))d\theta +M^q_u(\hat\xi^q_{0,x}(\theta))d \tilde w(\theta ),\quad \hat\xi_{0,x}(0)=x,
$$
$$   d\hat\beta^q(\theta)=\tilde c^q_u(\hat\xi^q(\theta))\hat\beta^q(\theta)d\theta+ C^q_u(\hat\xi^q(\theta)\hat\beta^q(\theta)dw(\theta),
\quad \hat\beta^q(0)=I
  $$
with coefficients $\tilde c^q,  C^q$ given by   (\ref{3.16})-- (\ref{3.18})   and
$$
 \pmatrix{u^q(t,x)\cr\nabla u^q(t,x)} =E\left[\pmatrix{\hat\beta^q_{11}(\theta)&0\cr \hat\beta^q_{21}(\theta)&\hat\zeta^q_{22}(\theta)} \pmatrix{u_0^q(\hat\xi^q_{0,x}(\theta))\cr v^q_{0}(\hat\xi^q_{0,x}(\theta))}\right]
$$
we will formulate conditions to ensure that the functions $u^q(t,x)$ defined by the last relation exist and  give the required generalized solution of the problem (\ref{1.3}),  (\ref{1.4}). This will be done in a forthcoming paper.

{\bf Acknowledgement} The financial support of the RNF Grant 17-11-01136 is gratefully acknowledged.

\end{document}